\documentclass[11pt]{amsart}

\usepackage{amsmath,amsthm}
\usepackage{amssymb}
\usepackage{latexsym}
\usepackage{euscript}
\usepackage{graphics}
\usepackage[frame,ps,dvips,matrix,arrow,curve,rotate]{xy}

\newcommand{\Sp}{{\rm Sp}}
\newcommand{\Spf}{{\rm Spf}}

\newcommand{\rig}{{\rm rig}}
\newcommand{\AN}{{\rm an}}
\newcommand{\an}{{^\AN}}

\newcommand{\Fr}{{\rm Fr}}

\newcommand{\U}{{\rm U}}

\newcommand{\ba}{\left( \begin{array}{cc}}
\newcommand{\ea}{\end{array} \right)}
\newcommand{\inc}{\hookrightarrow}
\newcommand{\liminv}{\varprojlim}
\newcommand{\ts}{\otimes} 
\newcommand{\ra}{\rightarrow} 
\newcommand{\lra}{\longrightarrow} 

\newcommand{\llra}[1]{\stackrel{#1}{\lra}} 
\newcommand{\eqra}{\llra{\sim}} 

\newcommand{\Fp}{{\mathbb{F}}_p}
\newcommand{\Q}{\mathbb{Q}}

\newcommand{\Z}{\mathbb{Z}}
\newcommand{\C}{\mathbb{C}}

\newcommand{\Qp}{\mathbb{Q}_p}

\newcommand{\cH}{\mathcal H}

 \newcommand{\cO}{\mathcal O}

\newcommand{\omg}{\underline{\omega}}

\newcommand{\be}{\begin{eqnarray}}

\newcommand{\ee}{\end{eqnarray}}

\newtheorem{thm}{Theorem}[section]
\newtheorem{prop}[thm]{Proposition}
\newtheorem{lemma}[thm]{Lemma}
\newtheorem{cor}[thm]{Corollary}

\theoremstyle{definition}
\newtheorem{defn}[thm]{Definition}

\theoremstyle{remark}
\newtheorem{rem}[thm]{Remark}

\newcommand{\cY}{\mathcal{Y}}

\newcommand{\cU}{\mathcal{U}}
\newcommand{\cV}{\mathcal{V}}

\newcommand{\pr}{{\rm pr}}

\newcommand{\cW}{\mathcal{W}}

\newcommand{\cE}{\mathcal{E}}

\newcommand{\cX}{\mathcal{X}}
\newcommand{\tX}{\tilde{X}}
\newcommand{\tZ}{\tilde{Z}}
\newcommand{\cZ}{\mathcal{Z}}

\newcommand{\tM}{\tilde{M}}

\newcommand{\Tate}{{\rm Tate}}

\begin{document}
\title{A Gluing Lemma And Overconvergent Modular Forms}

\author{Payman L Kassaei}

\address{Department of Mathematics\\
McGill University\\
Montreal, QC, H3A 2K6 \\ Canada}

\email{kassaei@math.mcgill.ca}

\begin{abstract}
We prove a gluing lemma for sections of line bundles on a rigid
analytic variety. We apply the lemma, in conjunction with a result
of Buzzard's \cite{Bu}, to give a proof of (a generalization of)
Coleman's theorem which states that overconvergent modular forms
of small slope are classical. The proof is ``geometric" in nature,
and is suitable for generalization to other Shimura varieties.
\end{abstract}

\maketitle

\section{Introduction}

Let~$p$ be a prime number,~$N\!\geq\! 5$ an integer prime to~$p$,
and~$m\!\geq\! 1$ and~$k$ integers. Let~$K$ be a finite extension
of~$\Qp$. Let~$X_1(Np^m)_K$ denote the modular curve of
level~$\Gamma_1(Np^m)$ over~$K$, and let~$\omg$ be the usual
invertible sheaf on $X_1(Np^m)_K$ which on the non-cuspidal locus
is the push-forward of the sheaf of invariant differentials of the
universal elliptic curve. Let $X_1(Np^m)_K^\an$ denote the rigid
analytic version of $X_1(Np^m)_K$, and denote the analytification
of $\omg$ again by $\omg$. An {\it overconvergent}~$p$-adic
modular form of level~$\Gamma_1(Np^m)$ and weight~$k$ defined
over~$K$ is a section of~$\omg^k$ on a rigid analytic subdomain
of~$X_1(Np^m)_K^\an$ which strictly contains the component of the
ordinary locus containing the cusp~$\infty$. The subspace of {\it
classical} modular forms consists of those sections of~$\omg^k$
which can be extended to~$X_1(Np^m)_K^\an$. By the rigid analytic
GAGA, any such section is the analytification of an honest modular
form (i.e., a section of $\omg^k$ on $X_1(Np^m)_K$). In
\cite{Co1,Co2}, Coleman proves the following.

\begin{thm}{\bf (Coleman \cite{Co1,Co2}).}\label{coleman} Let~$f$ be an
overconvergent~$\U_p$-eigenform of weight~$k$ with
eigenvalue~$a_p$. If the~$p$-adic valuation of~$a_p$ is less
than~$k-1$, then~$f$ is classical.
\end{thm}

This so-called ``control theorem" is crucial in many applications
of the theory to modular forms. In many problems, one can reduce
the ``rigidity" of the situation by working in Banach spaces of
overconvergent forms, and yet in the end get results about
classical modular forms, by invoking the control theorem. For
example, it is easier to construct~$p$-adic families of
overconvergent eigenforms, and then identify classical members of
these families. One is interested in such control theorems for
automorphic forms over other Shimura varieties. For instance, such
a result is still missing for (non-ordinary) overconvergent
Hilbert modular forms (See \cite{KL}). Coleman's method of proving
Theorem \ref{coleman}, which is based on a cohomological
interpretation of the space of overconvergent forms (modulo the
image of~$\theta^{k-1}$), appears difficult to be carried out in
more general situations (though it has the advantage that it gives
some information about the case of critical slope.)

The purpose of this paper is to present a more intrinsic proof of
(a generalization of) Theorem \ref{coleman} which could serve as a
model for extension to other Shimura varieties. In keeping with
this goal, we have tried to use the specifics of the situation as
little as possible. For instance, we use~$q$-expansions only to
show that modular forms are analytic at cusps.

In an upcoming article, we used this method to prove a control
theorem for overconvergent modular forms over other Shimura curves
(see \cite{Kas1,Kas2} for the general theory of overconvergent
modular forms over these Shimura curves). We expect this method to
apply without much trouble in some higher dimensional settings
(for instance, over unitary Shimura varieties which arise in the
work of Harris and Taylor \cite{HT}). The applicability of this
method in the case of Hilbert modular varieties, however, seems to
be contingent upon a better understanding (than at the moment) of
the canonical subgroups of Hilbert-Blumenthal Abelian varieties.
It also appears that our method can be applied in the context of
the theory of overconvergent modular forms of half-integral weight
developed by Nick Ramsey.

We were inspired by the work of Buzzard and Taylor \cite{BT} to
pursue the strategy of analytic continuation of modular forms. The
subsequent paper of Buzzard \cite{Bu} gave us the right starting
point. He proves that an overconvergent~$\U=\U_p$ eigenform $f$
with non-zero eigenvalue~$a=a_p$ can be extended over the
supersingular locus to a ``big" admissible open subset of the
modular curve. Under the assumption $v(a)<k-1$, we prove the
classicality of~$f$ by showing that~$f$ can be extended further,
over the complement of this admissible open. We prove the
following generalization of Theorem \ref{coleman}.
\begin{thm}{\label {classicalintro}} Let~$f$ be an overconvergent
modular form of weight~$k$ and level~$\Gamma_1(Np^m)$ defined over
$K$. Let~$R(x) \in K[x]$ be a polynomial all roots of which in
$\C_p$ have~$p$-adic valuation less than~$k-1$. If~$R(\U)f$ is
classical, then so is~$f$.
\end{thm}
Let us explain the steps involved in the proof of Theorem
\ref{classicalintro} in some detail. For the purpose of
presentation, we will take~$m=1$,~$R(x)=x-a$,~$v(a)\!<\!k-1$, and
we assume~$R(\U)f=0$, i.e.~$\U f=af$. Let~$\cZ^\infty(Np)$ (resp.,
$\cZ^0(Np)$) denote the connected component of the ordinary part
of $X_1(Np)_K^\an$ which contains the cusp~$\infty$ (resp.,~$0$).
By Buzzard's work \cite{Bu} one can extend~$f$ to~$\cU_1(Np)$
which is the rigid analytic part of~$X_1(Np)_K^\an$ whose
non-cuspidal points consist of all~$(E,i,P)$ where $E$ is an
elliptic curve, $i$ is a level $\Gamma_1(N)$-structure, and $P$ is
a point of order $p$ on $E$, and either~$E$ has supersingular
reduction, or~$E$ has ordinary reduction and~$P$ generates the
canonical subgroup of~$E$ (or equivalently~$(E,i,P) \in
\cZ^\infty(Np)$).

To show that~$f$ is classical, we have to show that~$f$ can be
further extended to the missing
part~$X_1(Np)_K^\an-\cU_1(Np)=\cZ^0(Np)$. This will be done in two
steps. The first step is to find a candidate for what the
extension of~$f$ to~$X_1(Np)_K^\an$ should be on~$\cZ^0(Np)$. The
second step is to show that this candidate indeed glues to (the
Buzzard extension of)~$f$. This step uses a general gluing lemma
that we prove in this paper, and involves a number of norm
estimations on modular forms using the theory of canonical
subgroups. The idea of the first part, however, is quite simple
(and fun!). Since in the paper this idea is presented in a rather
implicit way, we will give a clearer explanation for the benefit
of the reader.

To get an idea of how one should define~$f$ on~$\cZ^0(Np)$, let us
assume that~$f$ is classical for the moment, and see how the
values of~$f$ on~$\cZ^0(Np)$ are related to its values on the
complement of~$\cZ^0(Np)$ (where we know~$f$). We remark that for
a modular form $f$ and a test object~$(E,i,P)$ we think
of~$f(E,i,P)$ as an element of~$H^0(E,\Omega_{E})^{\ts k}$, {\`a}
la Katz. Assume~$(E,i,P)$ is in~$\cZ^0(Np)$. Since~$\U f= af$, we
can write
\begin{eqnarray}{\label{Uf=af}}
f(E,i,P)=(1/ap)\sum_{P \not\in C_1}
(\pr^*)^kf(E/C_1,\bar{i},\bar{P})
\end{eqnarray}where the sum is over the cyclic subgroups~$C_1$ of
order~$p=p^1$ which do not contain~$P$, and~$\pr:E \ra E/C_1$ is
the natural projection. By ``$\pr^*$'' we denote the pull-back of
one-forms under ``$\pr$''. One can show that all but one of the
test objects appearing on the right hand side of the above formula
belong to~$\cZ^\infty(Np)$. The exceptional term corresponds
to~$C=H_1$, the first canonical subgroup of~$E$. Applying the
above formula to~$(E/H_1, \bar{i},\bar{P})$ we get
\begin{eqnarray}
f(E/H_1,\bar{i},\bar{P})=(1/ap)\sum_{P \not\in C_2, H_1 \subset
C_2} (\pr^*)^kf(E/C_2,\bar{i},\bar{P})
\end{eqnarray}where the sum is over the cyclic subgroups~$C_2$ of
order~$p^2$ which contain~$H_1$, but not~$P$. Combining the two
equations we get
\begin{eqnarray}
f(E,i,P)=(1/ap)\sum_{P \not\in C_1,C_1 \neq H_1} (\pr^*)^k
f(E/C_1,\bar{i},\bar{P})+ (1/ap)^2\sum_{P \not\in C_2, H_1 \subset
C_2} (\pr^*)^k f(E/C_2,\bar{i},\bar{P})\nonumber.
\end{eqnarray}
We will repeat this process ad infinitum. At the~$n$-th step, we
separate the term corresponding to the quotient of~$E$ by~$H_n$,
the~$n$-th canonical subgroup, and rewrite the term via Equation
(\ref{Uf=af}). At the~$n$-th step the error term is~$(1/ap)^n
(\pr^*)^k f(E/H_n,\bar{i},\bar{P})$. Since on the (good-reduction
locus of the) ordinary part ~$H_n$ reduces to the kernel
of~$\Fr_{p}^n$ modulo~$p$, one can show that the error term is
divisible by~$(1/ap)^n p^{nk}=(p^{k-1}/a)^n$. Since~$v(a)<k-1$, we
see that the error term goes to zero as~$n$ goes to infinity. The
same estimate can be made on the locus of bad reduction as well
(see Corollary \ref{ordbnd}). The result is the following.

 \begin{thm}\label{formula}
Let~$f$ be a classical modular form of weight~$k$ and level
$\Gamma_1(Np)$ defined over~$K$. Assume that~$f$ satisfies~$\U
f=af$ with~$v(a)<k-1$. Then for~$(E,i,P) \in \cZ^0(Np)$ we have
\begin{eqnarray}
f(E,i,P)= \sum\limits_{n=1}^{\infty} (1/ap)^n {\big (}\sum_{C_n}
(\pr^*)^k f(E/C_n,\bar{i},\bar{P}) {\big )}
\end{eqnarray}where~$C_n$ runs through all the cyclic subgroups of~$E$
of order~$p^n$ which contain~$H_{n-1}$, are different from~$H_n$,
and do not contain~$P$.
\end{thm}
Now, let us go back to the assumption that~$f$ is overconvergent
of slope less than~$k-1$. One can show that all the test objects
appearing on the right hand side of the above formula belong to
$\cZ^\infty(Np)$, and hence,~$f$ is already defined for them. This
suggests that the extension of~$f$ to~$\cZ^0(Np)$, denoted by~$g$,
shall be defined via the above series. However this definition
will only work on~$ \cZ^0(Np)$ where~$H_n$ exists for all~$n$. The
rest is to prove that~$g$ which is defined on~$\cZ^0(Np)$, will
glue to~$f$ which is defined on the complement of~$\cZ^0(Np)$. And
that will follow from our gluing lemma, as we shall see.

Here goes the structure of this paper. In \S \ref{thelemma} we
prove the gluing lemma. In \S \ref{rev} we review some background
on the theory of modular curves and overconvergent modular forms,
and present (a slight generalization of) Buzzard's analytic
continuation results. In \S \ref{classicality} we prove Theorem
\ref{classicalintro}: we construct $g$ described above, and show
that the gluing lemma can be applied to glue $f$ and $g$ to
produce a classical modular form. To show that the lemma applies,
we obtain a variety of norm estimations using the theory of
canonical subgroups, and in particular, we prove the boundedness
of the Buzzard extension of~$f$ on the wide open
space~$\cU_1(Np^m)$.

{\bf Acknowledgment.} We are grateful to Kevin Buzzard and Richard
Taylor for their paper \cite{BT} which inspired us to pursue the
strategy of analytic continuation. More thanks are due to Buzzard
for his paper \cite{Bu}, the influence of which on this work is
clear. We also thank Robert Coleman and Ofer Gabber for
interesting discussions. Finally, we thank the referee for a very
close reading of this article, and the many useful suggestions and
comments that helped improve the presentation of this manuscript.

\section{The Gluing Lemma}\label{thelemma}

In this section we prove a gluing lemma for sections of line
bundles on rigid analytic varieties. Let~$K$ be a finite extension
of~$\Qp$ with the ring of integers~$\cO_K$. Let~$v$ denote the
valuation on~$K$ normalized to satisfy~$v(p)=1$. For~$c \in K$
define~$|c|=p^{-v(c)}$. Let~$X$ be a reduced flat scheme of finite
type over~$\cO_K$. Let~$X_K$ denote~$X \ts_{\cO_K} K$. Denote the
completion of~$X$ along its special fibre by~$\tX$. Note that
$\tX$ is flat and topologically finitely generated over~$\cO_K$
(i.e., an admissible formal scheme over~$\cO_K$). Raynaud's
functor associates to~$\tX$, its generic fibre~$\tX_\rig$, which
is a quasi-compact and quasi-separated rigid analytic space
over~$K$. One refers to~$\tX$ as a formal model for $\tX_\rig$.
See \cite{BL} for details on Raynaud's construction. There is also
the analytification functor which to~$X_K$ associates a rigid
analytic space~$X_K^\an$. In general, there is an open
immersion~$\tX_{\rig} \inc X_K^{\an}$ which is an isomorphism
if~$X$ is proper over~$\cO_K$.

Let~$M$ be an invertible sheaf on~$X$. This induces invertible
sheaves~$M_K$,~$\tilde{M}$,~$\tilde{M}_{\rig}$,~$M_K^\an$
respectively on~$X_K$,~$\tX$,~$\tX_\rig$,~$X_K^\an$. Let us fix a
finite trivialization of $\tilde{M}$ on $\tilde{X}$
\[
\{(\tilde{U}_i=\Spf(A_i),\tilde{\sigma}_i:\tM_{|_{\tilde{U_i}}}\eqra
\cO_{\tilde{U}_i}) \}_{i\in I}.
\]
This induces a trivialization for~$\tM_\rig$ on~$\tX_\rig$ which
we denote by~$\{ (\cU_i:=(\tilde{U}_i)_\rig ,\sigma_i)\}_{i\in
I}$.

\begin{defn}{\label{norm}} Fix the formal scheme $\tX$ and the
sheaf $\tM$ on it. Let~$x \in \tX_\rig$ be a point yielding a map
$x\colon\Sp(L)\ra \tX_\rig$, where $L$ is the residue field of
$x$. We first define a norm $|\ |_x$ on $H^0(\Sp(L),x^*\tM_\rig)$.
Denote the formal lifting of $x$ to the formal model by
$\tilde{x}:\Spf(\cO_L) \ra \tX$, where $\cO_L$ is the ring of
integers in $L$. Then
\[
H^0(\Sp(L),x^*\tM_\rig)=H^0(\Spf(\cO_L),\tilde{x}^*\tM)\ts_{\cO_L}
L
\]
and we define $|\ |_x$ via identifying
$H^0(\Spf(\cO_L),\tilde{x}^*\tM)$ with $\cO_L$. Clearly,  the
definition is independent of the identification. Now, consider an
admissible open subset~$\cU \subset \tX_\rig$, and let~$f \in
H^0(\cU,\tM_\rig)$ and $x \in \cU$. We define
\[
|f(x)|:=|x^*f|_x.
\]
We also define the norm of $f$ over $\cU$ (possibly infinite) to
be
\[
|f|_{_\cU}:=\sup\{|f(x)|: x \in \cU\}.
\]
If $x \in \cU \cap \cU_i$, we can use the identification
$H^0(\Spf(\cO_L),\tilde{x}^*\tM)\cong\cO_L$ induced by
$\tilde{\sigma_i}$ to calculate $|f(x)|$. This shows that
\begin{eqnarray}\label{normsup}
|f|_{_\cU}=\max_{i \in I} \{ |\sigma_i(f_{|_{\cU \cap
\cU_i}})|_{sup}\}
\end{eqnarray}
where the~$|\ |_{sup}$ on the right is the usual supremum
norm of functions.
\end{defn}

\

\begin{lemma}{\label{Banach}}
If~$\ \cU$  is an affinoid subdomain of~$\tX_\rig$, then~$|\
|_{_\cU}$ is a norm on~$H^0(\cU,\tM_\rig)$ which makes it into
a~$K$-Banach module.
\end{lemma}

\begin{proof}
First we show that~$|\ |_{_\cU}$ is finite. Since~$\tX_\rig$ is
quasi-separated,  each~$\cU \cap \cU_i$ is quasi-compact and
hence, it can be covered by finitely many affinoids. Therefore,
the supremum norm of functions is finite on each ~$\cU \cap
\cU_i$. Now, the result follows from Equation (\ref{normsup})
since $I$ is finite.

In general~$|\ |_{sup}$ gives a complete and separated norm on
functions on a reduced affinoid. The fact that~$\tX_\rig$ is
reduced along with Equation (\ref{normsup}) implies that~$|\
|_{_\cU}$ is separated.

For the completeness, note that a Cauchy sequence consisting of
elements in~$H^0(\cU,\tM_\rig)$ produces, for each~$i \in I$, a
Cauchy sequence of functions on the quasi-compact~$\cU \cap
\cU_i$.  Each of these Cauchy sequences of functions converges,
and the limits can be used, via the trivialization in question, to
produce a section of~$\tM_\rig$ on~$\cU$ which lies in the limit
of the original Cauchy sequence.
\end{proof}

Now, we prove our gluing lemma.

\begin{lemma}{\bf (The Gluing Lemma).}{\label {gluing}}
Let the notation be as above. Let~$\cX \subset \tX_\rig$ be a
smooth affinoid subdomain. Assume that~$\cX$ is a disjoint union
of two admissible opens~$\cX=\cY \cup \cZ$, where $\cZ$ is an
affinoid. Assume we are given affinoid subdomains of $\cX$ denoted
by~$\cZ_n$ for~$n \geq 1$ with
\[
\cZ \subset \cZ_1 \subset \cZ_2 \subset \cdots
\]and such that~$\{\cY,\cZ_n\}$ is an admissible cover of~$\cX$ for each
$n$. Assume that we are given two sections
\[
f\in H^0(\cY,\tM_\rig) \ \ \ \ \ \ \ \ \ g\in H^0(\cZ,\tM_\rig)
\]
and for each~$n \geq 1$, a section~$F_n \in H^0(\cZ_n,\tM_\rig)$
such that, as~$n \ra \infty$, we have
\[
|F_n-f|_{_{\cY \cap \cZ_n}} \ra 0\ \ \ {\rm and}\ \ \
|F_n-g|_{_\cZ} \ra 0.
\]
Then~$f$ and~$g$ glue together to give a global section
of~$\tM_\rig$ on~$\cX$. In other words, there is a section
of~$\tM_\rig$ on~$\cX$, which restricts to~$f$ on~$\cY$, and
restricts to~$g$ on~$\cZ$.
\end{lemma}

\begin{proof}

Since gluing $f$ and $g$ is a local issue, and in view of Equation
(\ref{normsup}), we can assume that $\tM_{\rig}$ restricted to
$\cX$ is the structure sheaf, and that for a section $f$ (over an
open of $\cX$), and for $x \in \cX$, the norm, $|f(x)|$, is the
same as the spectral norm for functions. From now on, we will
think of $f,g$, and all $F_n$'s as analytic functions, and of all
norms as supremum norms.

Let~$\check{\cO}_\cX$ denote the subsheaf of~$\cO_\cX$ whose
sections over an admissible open are those analytic functions with
supremum norm less than $1$.  It is not difficult to see that the
assumptions of the lemma imply that $|f|_{_\cY}$ and $|g|_{_\cZ}$
are finite, and $|F_n|_{_{\cZ_n}}$ are bounded independently of
$n$, and hence, by rescaling, we can assume that~$f$,~$g$, and
all~$F_n$'s are sections of~$\check{\cO}_\cX$. We can also assume
that
\begin{eqnarray} {\label{approximate}}|F_n-f|_{_{\cY \cap \cZ_n}} < (1/p)^n \ \ \ {\rm and}\ \ \
|F_n-g|_{_{\cZ}} < (1/p)^n.
\end{eqnarray}
Let~$\underline{\check{\cO}_\cX/p^n \check{\cO}_\cX}$ denote the
quotient sheaf. Then the reduction of~$f$ (resp.,~$F_n$)
modulo~$p^n$ is a section of~$\underline{\check{\cO}_\cX/p^n
\check{\cO}_\cX}$ on~$\cY$ (resp.,~$\cZ_n$) and by Equation
(\ref{approximate}) they agree over~$\cY \cap \cZ_n$. This implies
that they glue together to give a section~$h_n$
of~$\underline{\check{\cO}_\cX/p^n \check{\cO}_\cX}$ over~$\cX$.
The conditions of the lemma imply that~$h_n$'s are compatible in
the sense that they give an element of the inverse limit
$$\liminv_n \underline{\check{\cO}_\cX/p^n \check{\cO}_\cX}(\cX).$$ Theorem 2
of \cite{Ba} tells us that, since $\cX$ is a smooth affinoid,
there exists a non-zero $c \in K$ with $|c|\leq 1$ such that
$cH^1(\cX,\check{\cO}_\cX)=0$. This implies that
\[
ch_n\in \check{\cO}_\cX(\cX)/p^n\check{\cO}_\cX(\cX) \subset
\underline{\check{\cO}_\cX/p^n \check{\cO}_\cX}(\cX), \] and hence
$\{ch_n\}$ defines an element of $\liminv_n
\check{\cO}_\cX(\cX)/p^n\check{\cO}_\cX(\cX)=\check{\cO}_\cX(\cX)$
(since $\check{\cO}_\cX(\cX)$ is $p$-adically complete). We define
$h$ to be the section of $\cO_\cX$ over $\cX$ obtained by dividing
the above section by $c$.

By definition of~$h$, we have~$ch_{|_\cY}-cf \in
p^n\check{\cO}_\cX(\cY)$ for all~$n\geq 1$, and hence, we find
that~$h$ restricts to~$f$ on~$\cY$. Similarly, we find
that~$ch_{|_\cZ}-cg=(ch_{|_\cZ}-{cF_n}_{|_\cZ})+({cF_n}_{|_\cZ}-cg)
\in p^n\check{\cO}_\cX(\cZ)$ for all~$n\geq 1$ which implies
that~$h$ restricts to~$g$ over~$\cZ$. We are done.
\end{proof}

\begin{rem} The smoothness hypothesis can be removed at the
expense of making the lemma more technically involved. Since in
the applications to  overconvergent modular forms the smoothness
can be afforded we have chosen to make this assumption. We are
grateful to Ofer Gabber for pointing to us Bartenwerfer's result
on the cohomology of the sheaf $\check{\cO}_\cX$.
\end{rem}

\section{Overconvergent modular forms}\label{rev}

We review some background on the theory of overconvergent modular
forms. For more details we refer the reader to consult \cite{KM},
\cite{Bu}.

\subsection{Modular curves and modular forms} Our main reference is \cite{KM}.
Let~$N\!>\!4$ be an integer and~$p$ be a prime number such that
$(p,N)=1$. By~$X_1(N)$ we denote the smooth and proper modular
curve over~$\Z[1/N]$ whose non-cuspidal part, denoted~$Y_1(N)$,
classifies pairs~$(E,i)$, where~$E$ is an elliptic curve over
a~$\Z[1/N]$-scheme, and~$i\colon\mu_N \ra E$ is an embedding of
the constant group scheme~$\mu_N$ in~$E$. For
any~$\Z[1/N]$-algebra~$R$, we let~$X_1(N)_R$ denote the base
extension of~$X_1(N)$ to~$R$. Denote by~$E_1(N)$ the universal
family of elliptic curves over~$Y_1(N)$. There is a well-known
locally free sheaf of rank one,~$\omg=\omg_{X_1(N)}$, whose
restriction to~$Y_1(N)$ is the push-forward
of~$\Omega^1_{E_1(N)/Y_1(N)}$ under the natural map from~$E_1(N)$
to~$Y_1(N)$. The space of modular forms of weight~$k \in \Z_{\geq
0}$ and level~$\Gamma_1(N)$ over~$R$, denoted
by~$M_k(\Gamma_1(N),R)$, is~$H^0(X_1(N)_R, \omg^k)$.

Let~$X_1(N,p)$ denote the flat and proper modular curve over
$\Z[1/N]$ whose non-cuspidal part,~$Y_1(N,p)$, classifies
triples~$(E,i,C)$, where~$(E,i)$ is as above, and~$C$ is a finite
flat subgroup of~$E$ of order~$p$. For any~$\Z[1/N]$-algebra, $R$,
we let~$X_1(N,p)_R$ denote the base extension of~$X_1(N,p)$
to~$R$. As in the above, there is a universal elliptic
curve~$E_1(N,p)$ over~$Y_1(N,p)$, and a locally free
sheaf~$\omg=\omg_{X_1(N,p)}$ over~$X_1(N,p)$ whose restriction
to~$Y_1(N,p)$ is the push-forward of
$\Omega^1_{E_1(N,p)/Y_1(N,p)}$. If we let~$\pi_1\colon X_1(N,p)
\ra X_1(N)$ denote the degeneracy map which (on the non-cuspidal
part) forgets the subgroup of order~$p$, then one has
\[
\omg_{X_1(N,p)}=\pi_1^*\omg_{X_1(N)}.
\]
If $p$ is invertible in $R$, we define the space of modular forms
of weight~$k \in \Z_{\geq 0}$ and level~$\Gamma_1(N) \cap
\Gamma_0(p)$ over~$R$, denoted by~$M_k(\Gamma_1(N) \cap
\Gamma_0(p),R)$, to be~$H^0(X_1(N,p)_R, \omg^k)$.

Assume~$m>0$. Let~$X_1(Np^m)$ denote the proper and flat modular
curve over~$\Z[1/N]$, the non-cuspidal part of which, denoted
$Y_1(Np^m)$, classifies triples~$(E,i,P)$, where~$(E,i)$ is as
above, and~$P$ is a point of exact order~$p^m$. For
any~$\Z[1/N]$-algebra~$R$, we let~$X_1(Np^m)_R$ denote the base
extension of~$X_1(Np^m)$ to~$R$. If~$p$ is invertible in~$R$,
then~$X_1(Np^m)_R$ is smooth. As in the case of~$X_1(N)$, there
 a universal family of elliptic curves~$E_1(Np^m)$
over~$Y_1(Np^m)$, and a locally free sheaf of rank
one,~$\omg=\omg_{X_1(Np^m)}$, whose restriction to~$Y_1(Np^m)$ is
the push-forward of~$\Omega^1_{E_1(Np^m)/Y_1(Np^m)}$. The pullback
of~$\omg_{X_1(N)}$ under the natural map from~$X_1(Np^m)$ to
$X_1(N)$ which (on the non-cuspidal part) forgets the point of
order~$p^m$ is~$\omg_{X_1(Np^m)}$. When $p$ is invertible in $R$,
the space of modular forms of weight~$k \in \Z_{\geq 0}$ and
level~$\Gamma_1(Np^m)$ over~$R$, denoted
by~$M_k(\Gamma_1(Np^m),R)$, is~$H^0(X_1(Np^m)_R, \omg^k)$.

\subsection{Overconvergent modular forms}

Let~$K/\Q_p$ be a finite extension of~$\Q_p$ with valuation~$v$
normalized so that~$v(p)=1$, and with the corresponding norm~$|\
|=(1/p)^v$. Let~$\cO_K$ denote the ring of integers, and $\kappa$
the residue field.

Let~$X$ denote any of the modular curves described above. There is
a rigid analytic space~$X_K^{\an}$ associated to the modular curve
$X_K$ over~$K$. We will denote the analytification of the sheaf
$\omg$ again by~$\omg$. Let~$\cU$ be an admissible open subset
of~$X_K^\an$, and let~$f \in H^0(\cU,\omg^k)$. If~$x \in \cU$ is a
point, we define~$|f(x)|$ as in Definition \ref{norm} using the
admissible formal scheme~$\tX_{\cO_K}$, which is the completion of
$X_{\cO_K}$ along its special fibre, and the induced invertible
sheaf~$\omg^k$ on $\tX_{\cO_K}$. When $p>3$, the Eisenstein series
$E_{p-1}$ is a modular form of weight~$p-1$ which lifts the Hasse
invariant from characteristic $p$. For $x \in X_K^\an$ we define
the ``measure of supersingularity" of $x$ to be $|E_{p-1}(x)|$. In
general, and to include the cases $p=2,3$, the measure of
supersingularity of $x$ can  essentially be defined via the norm
of a parameter on the completion of $X_1(N)_{\cO_K}$ at the
reduction of the ``point" which induces the prime-to-$p$ part of
$x$. This coincides with $|E_{p-1}(x)|$ when $p>3$ and
$|E_{p-1}(x)|<1$. For details we refer the reader to \S 3 of
\cite{Bu} (where things are defined using valuation rather than
norm). In this paper, we always use $|E_{p-1}(x)|$ for the measure
of supersingularity, and keep in mind its generalized sense when
$p=2,3$.

 We will define various rigid
analytic subdomains of~$X_K^\an$. Let~$1/p < r \in |\C_p|$.
Let~$X_1(N)^{\geq r}_K$ be the affinoid subdomain
of~$X_1(N)_K^\an$ whose points are given by
\begin{eqnarray*}
\{x \in X_1(N)_K^\an\colon |E_{p-1}(x)| \geq r\}.
\end{eqnarray*}
The space of overconvergent modular forms of weight~$k$ and
level~$\Gamma_1(N)$ over~$K$ is
\begin{eqnarray*}
M_k^\dagger(N,K)=\varinjlim_{p^{-p/(p+1)}<r \ra 1^-}
H^0(X_1(N)_K^{\geq r}, \omg^k).
\end{eqnarray*}

If~$(E,i)/L$ is a point on~$X_1(N)_K^\an$ such
that~$|E_{p-1}(E,i)| > p^{-p/(p+1)}$ (i.e., if~$(E,i)$ is {\it not
too supersingular}), then~$E[p]$ has a canonical
subgroup~$H_1=H_1(E)$ of order~$p$, which is defined over~$\cO_L$
if~$E$ has good reduction, and serves as a canonical lifting of
the kernel of Frobenius. With more restriction on the
supersingularity of~$(E,i)$ we can define canonical subgroups of
higher order. For example, if~$|E_{p-1}(E,i)| > p^{-1/(p+1)}$, one
can show that~$|E_{p-1}(E/H_1(E),\bar{i})| > p^{-p/(p+1)}$, and
hence~$E/H_1(E)$ has a canonical subgroup of order~$p$. The
inverse image of this subgroup in~$E$ is a cyclic subgroup of
order~$p^2$ of~$E$, denoted by~$H_2(E)$, which contains~$H_1(E)$.
Similarly if~$|E_{p-1}(E,i)| > p ^{-p^{2-n}/(p+1)}$, then~$E$ has
a cyclic canonical subgroup of order~$p^n$, denoted by~$H_n(E)$,
and one has
\[
H_1(E) \subset H_2(E) \subset ... \subset H_n(E).
\]
For a detailed analysis of canonical subgroups and the history of
the subject see \S 3 of \cite{Bu}.

Let~$r> p^{-p/(p+1)}$. The map which (on the non-cuspidal part)
sends~$(E,i)$ to~$(E,i,H_1)$ provides a section for the forgetful
morphism~$\pi_1^\an \colon X_1(N,p)_K^\an \ra X_1(N)_K^\an$
over~$X_1(N)_K^{\geq r}$, and hence, gives an isomorphism
between~$X_1(N)_K^{\geq r}$ and its image in~$X_1(N,p)_K^\an$. The
image is denoted by~$X_1(N,p)_K^{\geq r}$, and is the connected
component of the cusp~$\infty$ in the affinoid subdomain
of~$X_1(N,p)_K$ defined by~$|E_{p-1}(x)| \geq r$. Its non-cuspidal
points consist of all~$(E,i,C)$ such that~$|E_{p-1}(E,i)| \geq r$
and~$C=H_1(E)$. This shows that there is an inclusion
\begin{eqnarray*}
H^0(X_1(N,p)_K^\an,\omg^k) \inc M_k^\dagger(N,K)
\end{eqnarray*}
and we call the modular forms in the image {\it classical} (over
$X_1(N,p)_K$). Indeed, since~$X_1(N,p)$ is proper over~$\Z[1/N]$,
we know that the analytification map
\[
H^0(X_1(N,p)_K,\omg^k)
\llra{\rm an} H^0(X_1(N,p)_K^\an, \omg^k)
\]
is an isomorphism.

We define these concepts for level~$\Gamma_1(Np^m)$ where~$m>0$.
Consider the map
\[
X_1(Np^m)_K \ra X_1(N,p)_K
\]
which on the non-cuspidal part sends~$(E,i,P)$ to~$(E/\langle
pP\rangle,\bar{i},\bar{P})$. Let~$r$ be an element of $|\C_p|$
such that $r>p^{-p^{2-m}/(p+1)}$. Define~$X_1(Np^m)_K^{\geq r}$ to
be the inverse image of $X_1(N,p)_K^{\geq r^{p^{m-1}}}$ under this
map. The non-cuspidal part of~$X_1(Np^m)_K^{\geq r}$ consists of
all~$(E,i,P)$ such that $|E_{p-1}(E,i)|\geq r$,
and~$H_m(E)=\langle P\rangle$. The space of overconvergent modular
forms of weight~$k$ and level $\Gamma_1(Np^m)$ over~$K$ is
\begin{eqnarray*}
M_k^\dagger(Np^m,K)=\varinjlim_{(1/p)^{p^{2-m}/(p+1)}<r \ra 1^-}
H^0(X_1(Np^m)_K^{\geq r}, \omg^k).
\end{eqnarray*}
There is, therefore, an inclusion
\[
H^0(X_1(Np^m)_K^\an,\omg^k) \inc M_k^\dagger(Np^m,K).
\]
An overconvergent modular form in~$M_k^\dagger(Np^m,K)$ is said to
be {\it classical} (over~$X_1(Np^m)_K$) if it is in the image of
the above inclusion, i.e., if it can be extended to a section
of~$\omg^k$ on~$X_1(Np^m)_K^\an$.

By the~$q$-expansion of an overconvergent modular form we mean
the~$q$-expansion at the cusp~$\infty$. There is a Hecke
operator~$\U_p=\U$ acting on~$M_k^\dagger(Np^m,K)$ which on
the~$q$-expansions has the effect
\begin{eqnarray*}
\U(\sum a_nq^n)=\sum a_{np}q^n.
\end{eqnarray*}
By a generalized eigenform of~$\U$ we mean an overconvergent
modular form~$f$, for which there is~$R(x) \in K[x]$ such that
$R(\U) f=0$. We say that~$f$ has slope~$\alpha \in \Q$, if all the
roots of~$R(x)$ in~$\C_p$ have valuation~$\alpha$. For instance,
the slope of a~$\U$-eigenform is the valuation of its eigenvalue.
One knows that the slope of a classical~$\U$-eigenform of level
$\Gamma_1(Np)$ and weight~$k$ is at most~$k-1$.

In this paper we think about modular forms as Katz does in
\cite{Ka}, albeit in a rigid analytic sense. In particular, if~$f$
is a modular form of weight~$k$ and level~$\Gamma$ defined over an
admissible open~$\cU$ of the corresponding modular curve,
and~$(E,\gamma)\in \cU$ is an elliptic curve with
level~$\Gamma$-structure, then~$f(E,\gamma) \in
H^0(E,\Omega_E)^{\ts k}$.

\subsection{Buzzard's analytic continuation results}
In this subsection, we will recall the analytic continuation
results obtained by Buzzard in his paper \cite{Bu}. Let
$\cZ^0(N,p)$ be the connected component of the cusp~$0$ in the
ordinary part of~$X_1(N,p)_K^\an$. Let~$\cU_1(N,p)$ be the
complement of~$\cZ^0(N,p)$. It is an admissible open
in~$X_1(N,p)_K^\an$ whose non-cuspidal points consist of~$(E,i,C)$
such that either~$E$ has supersingular reduction, or~$E$ has
ordinary reduction and~$C=H_1(E)$. This is denoted by~$\cW_0(p)$
in \cite{Bu}.

Assume~$m >0$. There is a map from~$X_1(Np^m)_K$ to~$X_1(N,p)_K$
which on the non-cuspidal part sends~$(E,i,P)$ to~$(E,i,\langle
p^{m-1}P\rangle)$. Let~$\cU_1(Np^m)$ denote the inverse image
of~$\cU_1(N,p)$ under this map. It is an admissible open
of~$X_1(Np^m)_K$ whose non-cuspidal points consist of
all~$(E,i,P)$ such that either~$E$ has supersingular reduction,
or~$E$ is of ordinary reduction and the subgroup generated by~$P$
contains~$H_1(E)$. In the case~$m=1$, Buzzard denotes~$\cU_1(Np)$
by~$\cW_1(p)$.

Buzzard shows that any overconvergent~$\U$-eigenform with non-zero
eigenvalue can be extended to~$\cU_1(Np^m)$. We will use a
generalization of this fact which can be proved by exactly the
same method.

\begin{thm}{\label{buz}}Let~$f$ be an overconvergent modular form of
level~$\Gamma_1(Np^m)$ over~$K$. Assume that for some~$R(x) \in
K[x]$ with~$R(0) \neq 0$ we can extend~$R(\U)f$ to~$\cU_1(Np^m)$ .
Then~$f$ can be extended to~$\cU_1(Np^m)$. The same statement is
valid in level~$\Gamma_1(N) \cap \Gamma_0(p)$.
\end{thm}

\begin{proof}
Let~$R(0)=-a$, and~$R_0(x)=R(x)+a$. Assume~$R(\U) f=F$ which is
defined over~$\cU_1(Np^m)$. The same method as in \cite{Bu} works:
only at the~$n$-th step instead of considering~$\U^n f/a^n$ as the
partial analytic continuation, one should take
\begin{eqnarray*}
R_0(\U)^n f /a^n - R_0(\U)^{n-1}F/a^{n} - ... -R_0(\U) F/a^2-F/a.
\end{eqnarray*}
One just needs to notice that~$Q(\U) F$ is defined over
$\cU_1(Np^m)$ for any~$Q(x)\in K[x]$ satisfying~$Q(0)=0$. A
similar argument works over~$X_1(N,p)_K$.
\end{proof}

\subsection{Some notation and a lemma} \label{notation}
Fix once and for all~$t \in |\C_p|$ with~$p^{-p/(p+1)}<t<1$. The
region in~$X_1(N,p)_K^\an$ where~$|E_{p-1}| \geq t$ has two
connected components: the connected component of the
cusp~$\infty$, and the connected component of the cusp~$0$ which
we denote by~$\cV_1(N,p)$. If~$(E,i,C)$ is a non-cuspidal point
of~$X_1(N,p)_K^\an$, then it lies in~$\cV_1(N,p)$
iff~$|E_{p-1}(E,i)| \geq t$ and~$C \neq H_1(E)$. This is an
affinoid subdomain of~$X_1(N,p)_K^\an$, and
\begin{eqnarray*}
\{\cV_1(N,p),\cU_1(N,p)\}
\end{eqnarray*}
is an admissible covering of~$X_1(N,p)_K^\an$. Assume~$m>0$. Let
\[
\phi\colon X_1(Np^m)_K^\an \ra X_1(N,p)_K^\an
\]
denote the map which sends~$(E,i,P)$ to~$(E,i,\langle
p^{m-1}P\rangle)$. Define~$\cV_1(Np^m)=\phi^{-1}\cV_1(N,p)$. Its
non-cuspidal points consist of all~$(E,i,P)$ such
that~$|E_{p-1}(E,i)| \geq t$ and the subgroup generated by~$P$
intersects~$H_1(E)$ trivially.
 Since~$\cU_1(Np^m)=\phi^{-1}\cU_1(N,p)$, we know that
\begin{eqnarray*}
\{\cV_1(Np^m),\cU_1(Np^m)\}
\end{eqnarray*}
is an admissible covering of~$X_1(Np^m)_K^\an$. We will use these
admissible coverings in \S \ref{classicality} when we deal with
classicality of overconvergent modular forms.

Let us fix the level~$\Gamma_1(Np^m)$ with~$m>0$. For simplicity,
we denote~$X_1(Np^m)_K^\an$ by~$X_K^\an$, and~$\cV_1(Np^m)$
by~$\cV$. For any closed, open, or half-open subinterval~$I$
of~$[0,1]$ (with endpoints in~$|\C_p|$), we define a corresponding
subdomain~$\cV I$ where~$|E_{p-1}|~$ falls within that interval.
For example~$\cV[r,s)$ denotes the part of~$\cV$ where~$r \leq
|E_{p-1}(x)| <s$. If $I$ is a closed interval, then $\cV I$ is an
affinoid. For any interval~$I$, let~$\cE I$ denote the universal
elliptic curve over~$\cV I-\{cusps\}$, and~$\cH_1 I$ denote the
canonical subgroup of~$\cE I$. By~$I^{p^n}$ we denote the interval
obtained by raising~$I$ to the power~$p^n$.
\begin{defn}{\label{pr}} Let~$r \in |\C_p|$ satisfy~$t \leq r \leq 1$. Let~$I$ be either~$[r,1]$ or~$[r,1)$.
Let
\begin{eqnarray}
 \tau\colon \cV I^{1/p} \ra \cV I \nonumber
\end{eqnarray}
be the map whose effect on the non-cuspidal points is given by
\[
\tau (E,i,P) = (E/H_1,\bar{i},\bar{P}).
\]
Here~$H_1$ is the canonical subgroup of~$E$ of
order~$p$,~$\bar{i}$ is the induced level~$\Gamma_1(N)$ structure
on~$E/H_1$, and~$\bar{P}$ is the image of~$P$ in~$E/H_1$. The
pullback of~$\omg_{|_{\cV I}}$ under~$\tau$ is an invertible
sheaf~$\omg^\prime$ on~$\cV I^{1/p}$ whose restriction to the
non-cuspidal part is ~$\omg_{(\cE I^{1/p}/\cH_1 I^{1/p})/(\cV
I^{1/p}-\{cusps\})}$. Around a cusp in $\cV$ the map $\tau$ is
given by $q \mapsto q^p$, and hence, the pullback under $\tau$ of
the canonical differential form on $\Tate(q)$, i.e., $dz/z$, is
the canonical differential form on $\Tate(q)/H_1(\Tate(q))\cong
\Tate(q^p)$, i.e., $dz^p/z^p$. This implies that $\omg^\prime$ is
the canonical extension (\`{a} la Katz-Mazur) of $\omg_{(\cE
I^{1/p}/\cH_1 I^{1/p})/(\cV I^{1/p}-\{cusps\})}$ to $\cV I^{1/p}$.
Therefore, there is a morphism of sheaves on~$\cV I^{1/p}$
\[
(\pr^*)^k:(\omg^\prime)^k \ra (\omg_{|_{\cV I^{1/p}}})^k
\]
which on the non-cuspidal part is induced by pulling back
differential forms under the natural projection~$\pr\colon \cE
I^{1/p} \ra \cE I^{1/p}/\cH_1 I^{1/p}$. Let~$f$ be a section
of~$\omg^k$ on~$\cV I$. We define~$f^\tau \in H^0(\cV I^{1/p},
\omg^k)$ by \[ f^\tau:=p^{-k}(\pr^*)^k\tau^*f.
\]
Hence one has
\begin{eqnarray}
f^\tau (E,i,P)=p^{-k} (\pr^*)^k f(E/H_1,\bar{i},\bar{P}) \nonumber
\end{eqnarray}where~$\pr\colon E \ra E/H_1$ is the natural projection.
All these can be defined over~$X_1(N,p)$ as well, in exactly the
same way.
\end{defn}

\

In the final part of this section, we prove a lemma which will be
used in the proof of Theorem \ref{classical}. Let~$Z_{\Z_p}$
denote the affine subscheme of $X_1(N)_{\Z_p}$ where $E_{p-1}$ is
invertible. When~$p=2$ or $3$, instead of $E_{p-1}$, we can use
$E_4$. The open non-cuspidal  subscheme of $Z_{\Z_p}$ represents
the functor which to every $\Z_p$-algebra,~$R$, associates an
elliptic curve over $R$ for which $E_{p-1}$ (or $E_4$ when
$p=2,3$) is invertible, with a level $\Gamma_1(N)$-structure. For
simplicity we  denote by $Z$ the base extension of $Z_{\Z_p}$ to
$\cO_K$. The formal completion of $Z$ along its special fibre,
$\tilde{Z}$, is a formal model for $\mathcal Z$, the locus of
ordinary reduction in $X_1(N)_K^\an$. Dividing by the canonical
subgroup on the non-cuspidal locus induces a map
\[
\sigma_n:Z \ts (\cO_K/p^n) \ra Z \ts (\cO_K/p^n)
\]
for each $n\geq1$, and by passing to the limit, a map
$\tilde{\sigma}: \tilde{Z} \ra \tilde{Z}$. This, in turn, induces
the Frobenius morphism $\sigma:\mathcal Z \ra \mathcal Z$. The
sheaf~$\omg$ on~$X_1(N)$ provides invertible sheaves on
$\tilde{Z}$ and $\mathcal Z$ which we again denote by $\omg$. Over
$\tZ$, we set $\omg^\prime:=\tilde{\sigma}^*\omg$; similarly, over
$\cZ$, we set $\omg^\prime:=\sigma^*\omg$. One can show that in
each case $\omg^\prime$ is the canonical extension to the cusps of
the sheaf of invariant differentials of the quotient of the
universal elliptic curve by its canonical subgroup.  There is,
therefore,  a map~$\pr^*: \omg^\prime \ra \omg$ on~$\tilde{Z}$
(and also $\mathcal Z$). Note that all the above constructions are
indeed defined over $\Z_p$ in the case of $\tilde{Z}$ (and serve
as formal models of the corresponding constructions in the case of
$\mathcal Z$). The natural forgetful map $X_1(Np^m)_K^\an \ra
X_1(N)_K^\an$ induces a map $\pi:\cV[1,1] \ra \cZ$ such that
$\pi^*\omg=\omg$ and $\pi^*\omg^\prime=\omg^\prime$. We also have
\begin{eqnarray}\label{compat1}
\sigma \circ \pi=\pi \circ \tau,
\end{eqnarray}
and
\begin{eqnarray}\label{compat2}
\pi^* \circ \pr^*=\pr^* \circ \pi^*
\end{eqnarray}
where the $\pr^*$ on the left is as above, and the one on the
right is as in Definition \ref{pr}.

\

\begin{lemma}{\label{tau}}
Let~$h \in H^0(\cV I,\omg^k)$, where~$I$ is either~$[r,1)$ or
$[r,1]$, and~$t\!\leq\!r\!\leq\!1$. For any point~$x \in \cV
I^{1/p}$ we have
\begin{eqnarray}
|h^\tau (x)|\leq |h(\tau x)||E_{p-1}(x)|^{-k}. \nonumber
\end{eqnarray}

\end{lemma}

\begin{proof} Let~$x\colon\Sp(L) \ra \cV I^{1/p}$ denote a point. The map
corresponding to the point~$\tau x$ is given by~$\tau
x\colon\Sp(L) \llra{x} \cV I^{1/p} \llra{\tau} \cV I$. It is
enough to show that
\[
|h(\tau x)|\leq 1 \Rightarrow |h^\tau(x)| \leq |E_{p-1}(x)|^{-k}.
\]
Let~$\tX$ denote the formal completion of~$X_1(Np^m)_{\cO_K}$
along its special fibre. Then $\tX_\rig=X_1(Np^m)_K^\an$.
Let~$\tilde{x}\colon\Spf(\cO_L) \ra \tX$ be the formal lifting
of~$x\colon\Sp(L) \llra{x} \cV I^{1/p} \ra X_1(Np^m)_K^\an$.
Similarly let~$\tilde{\tau x}$ denote the formal lifting of~$\tau
x$. We have
\[
x^*h^\tau=x^*(p^{-k}(\pr^*)^k\tau^*h)=p^{-k}(\pr_x^*)^k(\tau x)^*h
\]
where $(\pr^*)^k:H^0(\cV I^{1/p},(\omg^\prime)^k) \ra H^0(\cV
I^{1/p},\omg^k)$ is as in Definition \ref{pr}, and
\[
(\pr_x^*)^k:H^0(\Sp(L),x^*(\omg^\prime)^k=(\tau x)^*\omg^k) \ra
H^0( \Sp(L),x^*\omg^k)
\]
is the specialization of~$(\pr^*)^k$ to~$x$, and
satisfies~$x^*(\pr^*)^k=(\pr_x^*)^k x^*$. We claim that it
suffices to prove
\begin{eqnarray}{\label{pr^*}}
\pr_x^*(H^0(\Spf(\cO_L),(\tilde{\tau x})^*\omg))\subseteq
(p/E_{p-1}(x))H^0(\Spf(\cO_L),\tilde{x}^*\omg))
\end{eqnarray}
where~$E_{p-1}(x)$ is any element of~$L$ of norm~$|E_{p-1}(x)|$.
The reason is that
\[
|h(\tau x)|\leq 1 \Rightarrow |(\tau x)^*h|_{\tau x}\leq 1
\Rightarrow (\tau x)^*h \in H^0(\Spf(\cO_L),(\tilde{\tau
x})^*\omg^k)
\]
which in conjunction with (\ref{pr^*}) would give
\[
x^*h^\tau=p^{-k}(\pr_x^*)^k (\tau x)^*h \in
E_{p-1}(x)^{-k}H^0(\Spf(\cO_L),\tilde{x}^*\omg^k))
\]
which means that $|h^\tau(x)|=|x^*h^\tau|_x\leq |E_{p-1}(x)|^{-k}$
as desired. In what follows we prove (\ref{pr^*}).

\

First assume that~$x$ corresponds to an elliptic curve (with level
structure)~$\underline{E}/L$  which is supersingular (or more
generally of good reduction). By the assumptions in the lemma,~$E$
has a canonical subgroup~$H=H_1$. The map~$\pr_x^*$ can be
interpreted as
\[
\pr_E^*:H^0(E/H,\Omega_{E/H}) \ra H^0(E,{\Omega}_{E})
\]
where~$\pr_E:E \ra E/H$ is the natural projection.
Let~$\mathbf{E}$ be the integral model of~$E$ over~$\cO_L$, and
$\mathbf{H}$ the canonical subgroup of~$\mathbf{E}$. Let
$\pr_{\mathbf E}:\mathbf{E}\ra \mathbf{E/H}$ denote the
projection. To prove (\ref{pr^*}) for~$x$, it is enough to show
that
\[
\pr_{\mathbf E}^*H^0({\mathbf{E/H}},\Omega_{\mathbf{E/H}})\subset
(p/ E_{p-1}(\underline{\mathbf{E}},\omega))H^0(\mathbf{E},
\Omega_{\mathbf{E}})
\]
where~$\omega$ denotes any non-vanishing one-form on~$\mathbf{E}$.
By Theorem 3.1 in \cite{Ka}, the natural projection ~$\pr_{\mathbf
E}:\mathbf{E} \ra \mathbf{E}/\mathbf{H}$ reduces,
modulo~$p/E_{p-1}(\underline{\mathbf{E}},\omega)$, to~$\Fr_p$.
Since pulling back via~$\Fr_p$ kills one-forms the claim follows
in this case.

Next, we deal with $x \in \cV[1,1]$, i.e., with the
case~$|E_{p-1}(x)|=1$. Indeed we have already proven the statement
for all such~$x$ of good reduction. Recall the
map~$\pi\colon\cV[1,1] \ra \cZ$ defined in the paragraph before
this lemma. Let~$y=\pi(x)\in \cZ$. Denote
by~$\tilde{y}:\Spf(\cO_L)\ra \tZ$ the formal lifting of~$y$ to the
formal model~$\tZ$ of $\cZ$. Then the formal lifting of~$\pi(\tau
x)=\sigma y$ (see Equation (\ref{compat1}))
is~$\tilde{\sigma}\tilde{y}$. Since $\omg=\pi^*\omg$,
$\omg^\prime=\pi^*\omg^\prime$, and by Equation (\ref{compat2})
(and by integral versions of such compatibilities) it is
straightforward to see that (\ref{pr^*}) for~$x$ is equivalent to
\[
\pr_{\tilde{y}}^*(H^0(\Spf(\cO_L),(\tilde{\sigma}\tilde{y})^*\omg=\tilde{y}^*\omg^\prime))
\subset p (H^0(\Spf(\cO_L),\tilde{y}^*\omg))
\]
where~$\omg$ and~$\omg^\prime$ are on~$\tZ$, and
$\pr_{\tilde{y}}^*$ is the specialization of~$\pr^*:\omg^\prime
\ra \omg$ via~$\tilde{y}$. Hence, it suffices to show that
$\pr^*:\omg^\prime \ra \omg$ on $\tZ$ reduces to the~$0$ morphism
modulo $p$. Indeed, since this map is  defined over $\Z_p$, it is
enough to prove this statement over $\Z_p$. But on the
non-cuspidal part of~$\tZ_{\Z_p} \ts \Fp$ (which is the ordinary
part of~$X_1(N) \ts {\Fp}$) this map reduces, modulo~$p$,
to~$\Fr_p^*:(\omg \ts \Fp)^{(p)} \ra (\omg \ts \Fp)$ which is
zero. Now, a morphism of invertible sheaves on an integral curve
which is zero away from a finite number of points has to be zero.
We are done.
\end{proof}

\begin{rem} One can prove that the inequality in Lemma \ref{tau}
is indeed an equality, though we don't need it for our arguments.
\end{rem}

\

\begin{cor}\label{ordbnd}
Let~$h \in H^0(\cV[1,1],\omg^k)$. For all~$n \geq 1$ we have
\[
|h^{\tau^n}|_{_{\cV[1,1] }} \leq |h|_{_{\cV[1,1] }} < \infty.
\]
\end{cor}

\begin{proof}
This follows from Lemma \ref{tau} with $r=1$, and Lemma
\ref{Banach}.
\end{proof}

\

\section{Classicality of overconvergent modular forms} \label{classicality}

Recall that we are in level~$\Gamma_1(Np^m)$ with~$N\!>\!4$,
$(p,N)=1$, and~$m\!>\!0$, and we denote~$X_1(Np^m)_K^\an$,
$\cV_1(Np^m)$, and~$\cU_1(Np^m)$ by~$X_K^\an$,~$\cV$, and~$\cU$
respectively. We prove the following theorem.

\begin{thm}{\label{classical}} Let~$f$ be an overconvergent modular form
of weight~$k$ and level~$\Gamma_1(Np^m)$ defined over~$K$.
Let~$R(x)$ in~$K[x]$ be a polynomial all roots of which in~$\C_p$
have valuation less than~$k-1$. If~$R(\U)f$ is classical
on~$X_1(Np^m)_K$, then so is~$f$.
\end{thm}
Coleman's theorem follows as a result.
\begin{cor}{\bf{(Coleman)}}
Let~$f$ be a generalized~$\U$-eigenform  of weight~$k$, level
$\Gamma_1(Np^m)$, and slope less than~$k-1$. Then~$f$ is
classical.
\end{cor}
\begin{proof} If~$f$ is a generalized eigenform of weight~$k$ and slope less
than~$k-1$, then for some polynomial~$R(x)$ as in Theorem
\ref{classical}~$R(\U)f$ equals~$0$ which is very well classical!
Hence~$f$ is classical.
\end{proof}

\subsection{Proof of the Theorem}
We will now prove Theorem \ref{classical}. First we prove the
statement for~$R(x)$ of degree~$1$. Let~$f$ be an overconvergent
modular form of weight~$k$ on~$X_K^\an$ such that for
some~$a\!\in\! K$ of valuation less than~$k-1$, we can extend~$\U
f-af$ to a section~$F$ of~$\omg^k$ on~$X_K^\an$. By Buzzard's
theorem \ref{buz}, one can extend~$f$ to~$\cU$. Since~$\{\cU,\cV
\}$ is an admissible covering of~$X_K^\an$, to prove the
classicality of $f$, it is enough to show that the restriction
of~$f$ to~$\cV \cap \cU = \cV[t,1)$ (which we still denote by~$f$)
can be extended to~$\cV=\cV[t,1]$. Let~$b=p^{k-1}/a$. By
definition of~$\tau$, we know that~$f-bf^\tau$ is a section
of~$\omg^k$ on~$\cV[t^{1/p},1)$.

\begin{prop}{\label {extend}}
The section~$f-bf^\tau \in H^0(\cV[t^{1/p},1),\omg^k)$ extends to
a section~$F_1$ of~$\omg^k$ on~$\cV[t^{1/p},1]$.
\end{prop}
\begin{proof}
 First we define~$F_1$ away from the cusps. Let~$(E,i,P)\in \cV[t^{1/p},1]-\{cusps\}$. Let~$C$ be any subgroup
of~$E$ of order~$p$ which is different from~$H_1$ and which
intersects~$\langle P\rangle$ trivially, i.e.,~$C \neq \langle
p^{m-1}P\rangle$. Then the image of~$p^{m-1}P$ in~$E/C$ generates
the canonical subgroup of~$E/C$. The reason is that in general
if~$E$ has a canonical subgroup, then for any non-canonical
subgroup~$C$ of~$E$,~$E/C$ has a canonical subgroup which
equals~$E[p]/C$. See Theorem 3.10.7.3 in \cite{Ka}.
Therefore~$(E/C, \bar{i},\bar{P}) \in \cU$, and hence~$f(E/C,
\bar{i},\bar{P})$ is defined. Define
\begin{eqnarray}
F_1(E,i,P):=(1/ap)\sum_{C} (\pr^*)^k f(E/C, \bar{i},\bar{P}) -
(1/a) F(E,i,P)
\end{eqnarray}
where~$C$ runs through the subgroups of~$E$ of order~$p$, which
are different from~$H_1$, and intersect~$\langle P\rangle$
trivially. Then on~$\cV[t^{1/p},1]-\{cusps\}$ we have
\begin{eqnarray}
F_1(E,i,P)\!&\!=\!&\!(1/ap)\sum (\pr^*)^k f(E/D, \bar{i},\bar{P})-(1/ap) (\pr^*)^kf(E/H_1,\bar{i},\bar{P})- (1/a) F(E,i,P)\nonumber \\
 &=&\!(1/a) \U f (E,i,P)- bf^\tau (E,i,P)-(1/a) F(E,i,P)\nonumber \\
 &=&\!f(E,i,P)-bf^\tau(E,i,P)\nonumber
\end{eqnarray}
where in the above~$D$ runs through all subgroups of~$E$ of order
$p$ which intersect~$\langle P\rangle$ trivially. This proves that
$F_1$ extends~$f-bf^\tau$ on~$\cV[t^{1/p},1]-\{cusps\}$. To end
the proof, we show that~$F_1$ extends to~$\cV[t^{1/p},1]$. We
calculate the~$q$-expansion of~$F_1$ around the cusps in $\cV$.
For this, we can assume that $K$ contains a primitive $p^m$-th
root of unity $\zeta_{p^m}$. The cusps in $\cV$ consist of all
isomorphism classes of $(\Tate(q),i,q^{A/p^m}\zeta_{p^m}^B)$,
where $i:\mu_N \ra \Tate(q)$ is an inclusion, and $A$ is not
divisible by $p$. We show the analyticity around the cusp $c$
corresponding to $A=1,B=0$. The other cases are similar. Let
$\zeta_p=(\zeta_{p^m})^{p^{m-1}}$.
\begin{eqnarray}
F_1(c)&=&(1/ap)\sum_{j=1}^{p-1} (\pr^*)^kf(\Tate(q)/\langle
q^{1/p}\zeta_{p}^j\rangle,\bar{i},{\bar q}^{1/p^m})-(1/a)F(c)\nonumber\\
&=&(1/ap)\sum_{j=1}^{p-1}f(\Tate(q^{1/p}\zeta_{p}^j),i,q^{1/p^m})-(1/a)F(c)\nonumber
\end{eqnarray}
But $(\Tate(q^{1/p}\zeta_{p}^j),i,q^{1/p^m})$ can be obtained from
$c^\prime_j\!=\!(\Tate(q),i,q^{1/p^{m-1}}\zeta_{p^m}^{-j})$ via a
base extension sending $q$ to $q^{1/p}\zeta_p^j$. For $j$ not
divisible by $p$, we have $c^\prime_j \in \cU$. The reason is that
we have $(q^{1/p^{m-1}}\zeta_{p^m}^{-j})^{p^{m-1}}=\zeta_p^{-j}$
which generates the canonical subgroup of $\Tate(q)$. Therefore,
$f$ is analytic at $c^\prime_j$. Furthermore $F$ is also analytic
at $c$. Hence $f$ can be extended to the cusp $c$.
\end{proof}

We now set the stage for the application of our gluing lemma. For
each~$n\!\geq\!1$ we define a section~$F_n \in
H^0(\cV[t^{1/p^n},1],\omg^k)$ by
\[
F_n=\sum_{i=0}^{n-1}b^iF_1^{\tau^i}.
\]
By Proposition \ref{extend}
we have ${F_1}_{|_{\cV[t^{1/p},1)}}=f-bf^\tau$, and it follows
that on $\cV[t^{1/p^n},1)$ we have
\[
{F_n}_{|_{\cV[t^{1/p^n},1)}}=f-b^nf^{\tau^n}.
\]
Since~$v(b)>0$, Corollary  \ref{ordbnd} and Lemma \ref{Banach}
imply that
\[
g:=\sum_{i\geq 0} b^iF_1^{\tau^i} \nonumber
\]
converges to a section of $\omg^k$ on $\cV[1,1]$. This is exactly
the infinite sum that was met in the Introduction (for the
case~$m=1$). For each~$n\geq 1$ we have
\begin{eqnarray*}
{F_n}_{|_{\cV[1,1]}}=g-b^ng^{\tau^n}.
\end{eqnarray*}

We now want to apply Lemma \ref{gluing} to show that~$f$ and~$g$
glue together to form a section~$\mathbf{f}$ of~$\omg^k$ on~$\cV$.
Then~$\mathbf{f}$ (defined on~$\cV$) will agree over~$\cV \cap
\cU$ with~$f$ (defined on~$\cU$), and we will get the desired
extension of~$f$ to the whole modular curve.

We will apply the gluing lemma
with~$X=X_1(Np^m)_{\cO_K}$,~$M=\omg^k$,~$\cX=\cV$ which is a
smooth affinoid subdomain of $\tX_\rig=X_K^\an$,
$\cY=\cV[t,1)$,$\cZ=\cV[1,1]$,~$\cZ_n=\cV[t^{1/p^n},1]$. Since we
have~$F_n-f=-b^nf^{\tau^n}$,~$F_n-g=-b^ng^{\tau^n}$, and since
$v(b)=k-1-v(a)>0$, to verify the conditions of the Lemma, it's
enough to show that~$\{|g^{\tau^n}|_{_{\cV[1,1]}}\}_n$ and
~$\{|f^{\tau^n}|_{_{\cV[t^{1/p^n},1)}}\}_n$ are uniformly bounded.
The first assertion follows from Corollary \ref{ordbnd}. In the
following we prove the second assertion.

\begin{lemma}{\label{bound}}
The section~$f$ is bounded on~$\cV[t,1)$.
\end{lemma}

\begin{proof} Since~$f-bf^\tau$ extends to ~$\cV[t^{1/p},1]$ which
is an affinoid,~$|f-bf^\tau|_{_{\cV[t^{1/p},1)}}$ is finite by
Lemma \ref{Banach}. Similarly,~$|f|_{_{\cV[t,t^{1/p}]}}$ is finite
since $\cV[t,t^{1/p}]$ is an affinoid. Let~$M_1$ be a common upper
bound. We show, by induction on~$n$, that~$f$ is bounded
by~$M_1t^{-k(1/p+...+1/p^n)}$
on~$\cV_n\colon=\cV[t^{1/p^{n}},t^{1/p^{n+1}}]$ for all~$n$.
Let~$x \in \cV_{n+1}$. Then~$\tau(x) \in \cV_n$ and by the
induction hypothesis and Lemma \ref{tau} we have
\begin{eqnarray}
|f^\tau(x)| &\leq& |f(\tau(x))| |E_{p-1}(x)| ^{-k} \nonumber \\
&\leq & M_1t^{-k(p^{-1}+...+p^{-n})}t^{-kp^{-(n+1)}}=M_1t^{-k(p^{-1}+...+p^{-n}+p^{-(n+1)})}\nonumber
\end{eqnarray}
So we can write
 \begin{eqnarray}|f|_{_{\cV_{n+1}}} &\leq& \max\{ |f-bf^\tau|_{_{\cV_{n+1}}} , |b{f^\tau}|_{_{\cV_{n+1}}}  \} \nonumber \\
&\leq& \max\{M_1, |{f^\tau}|_{_{\cV_{n+1}}} \} \nonumber \\
&\leq& M_1t^{-k(p^{-1}+...+p^{-n}+p^{-(n+1)})}.\nonumber
\end{eqnarray}
Now it is clear that~$f$ on~$\cV[t,1)= \bigcup_{n\geq 1} \cV_n$ is
bounded by
\begin{eqnarray}
\sup_n \{ M_1t^{-k(p^{-1}+...+p^{-n}+p^{-(n+1)})} \}=
M_1t^{-k/(p-1)}=:M_2.\nonumber
\end{eqnarray}
\end{proof}
We are now able to prove the desired uniform boundedness.

\begin{lemma} There is an~$M>0$ such that for all~$n \geq 0$ we have
\[
|f^{\tau^n}|_{_{\cV[t^{1/p^n},1)}} \leq M.
\]
\end{lemma}
\begin{proof} Let~$x \in \cV[t^{1/p^n},1)$ .  By Lemma \ref{tau}, we have
\[
|f^{\tau^n}(x)| \leq
|f^{\tau^{n-1}}(\tau(x))| |E_{p-1}(x)| ^{-k}.
\]
Inductively, we find
\[
|f^{\tau^n}(x)| \leq |f(\tau^n(x))|
{\big (} \prod_{j=0}^{n-1} |E_{p-1}(\tau^{j}(x))| {\big )} ^{-k}.
\]
Now ~$x \in \cV [t^{1/p^n},1)$ satisfies~$|E_{p-1}(x)| \geq
t^{1/p^n}$. Theorem 3.10.7(2) in \cite{Ka} that we have
implies~$|E_{p-1}(\tau^{j}(x))| =|E_{p-1}(x)|^{p^j}$. Also, since
$\tau^n(x)\in \cV[t,1)$, by Lemma \ref{bound} we find
that~$|f(\tau^n(x))|\leq M_2$. Putting all this together, we have
\begin{eqnarray}
|f^{\tau^n}(x)| \leq M_2|E_{p-1}(x)|^{-k(p^n-1)/(p-1)}\leq M_2 t^{
-k(p^n-1)/p^n(p-1)} \leq M_2 t^{ -k/(p-1)}=:M.\nonumber
\end{eqnarray}

\end{proof}

This proves Theorem \ref{classical} for~$R(x)$ of degree one. Here
is how we deal with the general case: It's enough to prove the
theorem over a finite extension of~$K$, and hence, we can
assume~$R(x)=(x-a_1)(x-a_2)...(x-a_l)$, such that~$v(a_j)<k-1$ for
all~$j$. We have
\[
(\U-a_1)(\U-a_2)...(\U-a_l)f=F.
\]
Let~$f_j=(\U-a_{j+1})(\U-a_{j+2})...(\U-a_l)f$. Then for~$f_1$ we
have~$(\U-a_1)f_1=F$ and hence by the above~$f_1$ is classical.
For~$f_2$ we have~$(\U-a_2)f_2=f_1$ and since~$f_1$ is classical,
we deduce that~$f_2$ is classical. Continuing this way we get
that~$f$ is classical. The proof of Theorem \ref{classical} is now
complete.

\end{document}